\documentclass[12pt]{amsart}
\usepackage{graphicx} 
\usepackage{amsfonts, amsmath, amssymb, amsthm}
\usepackage{xcolor}
\usepackage{mathrsfs}
\usepackage{tcolorbox}
\usepackage{hyperref}
\usepackage{stix}
\usepackage{dirtytalk}
\usepackage[headheight=120pt,top=1.4in, bottom=1.4in, left=1.0in, right=1.0in]{geometry}
\title[Characterization of $H$- Valued Fourier Transform]{Characterization of Integral Transforms for Hilbert Space-Valued Functions}
\author[Achar, and Manna ]{Satyajyoti Achar, and Ramesh Manna }

\address{School of Mathematical Sciences, National Institute of Science Education and
Research, Bhubaneswar, An OCC of Homi Bhabha National Institute, Jatni 752050, India.}
\address{Homi Bhabha National Institute, Training School Complex, Anushakti Nagar, Mumbai 400094, India}
\email{satyajyoti.achar@niser.ac.in, rameshmanna@niser.ac.in}


\date{\today}

\date{\today}
\begin{document}
\newcommand{\R}{\mathbb{R}}
\newcommand{\Hsn}{\mathbb{H}}
\newcommand{\Sl}{\mathcal{L}}
\newcommand{\A}{\mathcal{A}}
\newcommand{\F}{\mathscr{F}}
\newcommand{\N}{\mathbb{N}}
\newcommand{\Hi}{\mathcal{H}}
\newcommand*{\bigchi}{\mbox{\Large$\chi$}}
\newcommand{\Hm}{\mathrm{H}}
\newcommand{\blp}{\big (}
\newcommand{\brp}{\big )}
\newcommand{\Z}{\mathbb{Z}}
\newcommand{\Za}{\mathscr{Z}}
\newcommand{\T}{\mathbb{T}}
\newcommand{\J}{\mathbf{J}}
\newcommand{\C}{\mathbb{C}}
\newcommand{\diff}{\mathrm{d}}
\newcommand{\Df}{{\mathcal{F}_T}}
\newcommand{\D}{\mathcal{D}}
\newcommand{\B}{\mathcal{B}}
\newcommand{\M}{\mathcal{M}}
\newcommand{\G}{\mathcal{G}}
\newcommand{\nc}{\circledast_A}
\newcommand{\tc}{\otimes_\theta^\lambda}
\newcommand{\im}{\mathrm{Im}}
\newcommand{\re}{\mathrm{Re}}
\newcommand{\ch}{\mathcal{C_\theta}}
\newcommand{\nch}{\mathcal{C_{-\theta}}}
\newcommand{\f}{\mathcal{F}}
\newcommand{\He}{\mathbb{H}^n}
\def\Msp{{\boldsymbol M}}
\newtheorem{theorem}{Theorem}[section]
\newtheorem{lemma}[theorem]{Lemma}

\theoremstyle{definition}
\newtheorem{definition}[theorem]{Definition}
\newtheorem{example}[theorem]{Example}
\newtheorem{OP}[theorem]{Open problem}
\newtheorem{AP}[theorem]{Appendix}
\newtheorem{corollary}[theorem]{Corollary}
\newtheorem{xca}[theorem]{Exercise}
\newtheorem{propn}[theorem]{Proposition}
\theoremstyle{remark}
\newtheorem{remark}[theorem]{Remark}

\numberwithin{equation}{section}

\begin{abstract}
The aim of this paper is to characterize the Fourier transform of functions taking values in a separable Hilbert space $\Hi$ and defined on locally compact abelian groups.
    Inspired by Krukowski’s characterization of integral transforms on locally compact abelian groups, we establish characterizations of the Fourier, Laplace, and cosine transforms for separable Hilbert space-valued functions on locally compact abelian groups. Furthermore, we obtain a characterization of the group Fourier transform on the Heisenberg group. 
\end{abstract}
\maketitle

\section{Introduction and Preliminaries}
The Fourier transform is a fundamental operator in harmonic analysis with far-reaching applications in mathematics, physics, and engineering, including partial differential equations, number theory, quantum mechanics, image processing, and infrared spectroscopy. A natural question is whether the Fourier transform can be characterized solely by its intrinsic properties. As observed by Krukowski \cite{krukowski1}, the answer is affirmative: the Fourier transform is uniquely determined by a suitable collection of algebraic and analytic axioms. This viewpoint, referred to as the \say{axiomatic} approach, seeks to characterize the Fourier transform by proving that any operator satisfying these prescribed properties necessarily agrees with the Fourier transform, up to a normalization constant.

One of the most fundamental properties of the Fourier transform is the convolution theorem, which states that the Fourier transform converts convolution into pointwise multiplication, and pointwise multiplication into convolution. This property has been the basis for several characterization results. Alesker \emph{et al.} \cite{alesker2008, alesker2009} initiated this line of research by characterizing the Fourier transform through its convolution property. Motivated by their work, Jaming \cite{Jaming} established analogous characterization theorems for the Fourier transform on the locally compact abelian groups $\mathbb{R}$, $\mathbb{T}$, $\mathbb{Z}$, and the finite cyclic group $\mathbb{Z}_n$. Subsequently, Biswas and Radha \cite{Hasan1} investigated the corresponding characterization problem for Hilbert space-valued functions. More recently, Krukowski \cite{krukowski1} developed a unified framework for characterizing the Fourier transform on arbitrary locally compact abelian groups solely through its convolution property.



In the same work, Krukowski \cite{krukowski1} further showed that the convolution property also characterizes other classical integral transforms. More precisely, he proved that the cosine and Laplace transforms are uniquely determined by the corresponding cosine and Laplace convolutions, respectively. 

The main objective of this paper is to investigate the characterization problem for integral transforms for separable Hilbert space-valued functions on locally compact abelian groups. Specifically, we characterize the Fourier, Laplace, and cosine transforms on spaces of separable Hilbert space-valued functions over locally compact abelian groups. 

We assume throughout that the Hilbert space $\Hi$ is separable. Our proof is based on the techniques developed by Krukowski \cite{krukowski2, krukowski1}.
A main difficulty in the Hilbert space-valued setting is that a single convolution is not sufficient. Indeed, for $f,g\in L^1(G,\Hi)$, we define (see \eqref{eq: vector conv def})
\[
f*g\in L^1(G,\mathbb C),
\]
and does not belong to $L^1(G,\Hi)$. We therefore introduce a second convolution operation, denoted by $\odot$, between scalar-valued and $\Hi$-valued functions (see \eqref{E3.4}). This ensures that
\[
(f*g)\odot h\in L^1(G,\Hi),\qquad f,g,h\in L^1(G,\Hi).
\]
Thus, unlike the classical scalar-valued setting, our framework involves the interplay of two distinct convolution operations. Consequently, the proof differs substantially from the classical case and requires additional arguments.




In \cite{Jaming}, Jaming posed a similar question for twisted convolution. By answering Jaming’s question affirmatively, Lavanya and Thangavelu characterized the Weyl transform. More precisely, they showed that any continuous $\ast$-homomorphism from $L^1(\mathbb{C}^d)$, equipped with twisted convolution as multiplication, into $B(L^2(\mathbb{R}^d))$ is, up to a natural equivalence, the Weyl transform. This result subsequently led to a characterization of the group Fourier transform on the Heisenberg group; see \cite{Lavanya_Heisen}. Inspired by this approach, Kumar and Sivananthan adapted the method of Lavanya and Thangavelu to characterize the Fourier transform on compact groups. More recently, Krukowski \cite{krukowski2} established a similar result using a different approach based on representation theory and the correspondence between representations of a compact group and its $L^1$-space. 

In the first part, we adapt Krukowski’s approach to obtain a characterization of the group Fourier transform on the Heisenberg group.
The Heisenberg group $\He= \C^n\times \R$ forms a noncommutative group under the operation 
$$(z,t)(w,s)=\left(z+w, t+s+\frac{1}{2}\mathfrak{I}\left(z.\overline{w}\right)\right)$$
The Haar measure on $\He$ is the Lebesgue measure $dz\, dt$ on $\C^n\times\R.$ Stone-von Neumann theorem guarantees that all the infinite-dimensional irreducible unitary representation of $\He$, acting on $L^2(\R^n)$, are parametrized by $\lambda\in \R^*$, and can be given by 
$$\pi_\lambda(z,t)\varphi(\xi)= e^{i\lambda t}e^{i\lambda\left(x.\xi+ \frac{1}{2}x.y\right)}\varphi\left(\xi+y\right),~~~~~~\xi\in \R^n, ~~~~\varphi\in L^2(\R^n).$$
For a detailed study of Heisenberg group, we refer to\cite{Thangavelu_Heisen_book},\cite{Geller}.
Now let $f:\He\to \mathbb{C}$ be an integrable function. Then the group Fourier transform of $f$ is defined as
$$\widehat{f}(\lambda)=\int_{\He}f\left(z,t\right)\pi_\lambda(z,t)\, dz\, dt, ~~~~ \lambda\in \R^*.$$
Let $\mathcal{B}\left(L^2(\R^n)\right),$ denote the space of all bounded linear maps from $L^2(\R^n)$ to $L^2(\R^n)$. Then we can show that $\widehat{f}(\lambda)\in \mathcal{B}\left(L^2(\R^n)\right),$ with $\|\widehat{f}(\lambda)\|_{op}\leq \|f\|_{L^1(\He)}.$
\vspace{5pt}

\noindent
For $f,g \in L^1(\He)$, the convolution $f\ast g$ is defined by
$$(f\ast g)(z,t)=\int_{\He}f((z,t)(-w,-s))g(w,s)\, dw\, ds, \quad (z,t)\in \He,$$
whenever the integral exits. Now it is easy to check the following properties of the group Fourier transform. 
\begin{itemize}
    \item [1.] $\left(\widehat{f^\ast}\right)(\lambda)=\left(\widehat{f}(\lambda)\right)^\ast$ for all $\lambda\in \R^\ast$, where $f^\ast(z,t)=\overline{f(-z,-t)}$ and $\left(\widehat{f}(\lambda)\right)^\ast$ is the adjoint of the operator in $\mathcal{B}(L^2(\R^n)).$
    \item [2.] $\left(f\ast g\right)^{\widehat{}}(\lambda)=\widehat{f}(\lambda)\widehat{g}(\lambda), \quad \lambda\in \R^\ast, \quad f,g \in L^1(\He).$
    \item [3.] $\left(R_{(z,t)}f\right)^{\widehat{}}(\lambda)= \widehat{f}(\lambda)\pi_\lambda(z,t)^\ast, \quad (z,t)\in \He$, where $R_{(z,t)}$ denotes the right translation operator, given by
    $$\left(R_{(z,t)}f\right)(w,s)= f\left(\left(w,s\right)\left( z,t\right)\right), \quad (w,s)\in \He.$$
\end{itemize}

Using this properties, Lavanya and Thangavelu characterized the group Fourier transform.
Now we state their result.

\begin{theorem}[\cite{Lavanya_Heisen}]  \label{T1.1}
    Let $T: L^1(\He)\to L^\infty(\R^\ast, \mathcal{B}\left(L^2(\R^n)\right), d\mu)$ be a nonzero continuous linear map satisfying 
    \begin{itemize}
        \item [(i)] $T(f^\ast)(\lambda)=Tf(\lambda)^\ast,~$ for all  $ \lambda\in \R^\ast, f\in L^1(\He),$
        \item [(ii)] $T(f\ast g)(\lambda)=(Tf)(\lambda)(Tg)(\lambda), ~$ $\lambda\in \R^\ast, ~ f,g \in L^1(\He),$ and
        \item [(iii)] $T(R_{(z,t)}f)(\lambda)= Tf(\lambda)\pi_\lambda(z,t)^\ast, ~ \lambda\in \R^\ast, f\in L^1(\He), t\in \R.$
    \end{itemize}
    
    Then $$\left(Tf\right)(\lambda)=\widehat{f}(\lambda),~ \lambda\in A, f\in L^1(\He),$$
    where $A:=\{\lambda\in \R^\ast: (Tf)(\lambda)\neq 0 \text{ for some } f\in L^1(\He)\}.$
\end{theorem}
Now we see the characterization results on the compact groups. For that first we define some preliminaries for the Fourier analysis on the compact groups. Let $G$ be a compact group, not necessarily abelian. Then any irreducible unitary representation of $G$ is always finite dimensional. Let $\widehat{G}$ denote the set of all unitary equivalence classes of irreducible unitary representation of $G$. It is called the unitary dual of $G$, given the discrete topology.

We shall denote by $l^\infty- \bigoplus_
{\alpha \in \Delta}X_\alpha$, the Banach space $\{(x_\alpha)\in \prod_{\alpha\in \Delta}X_\alpha: \sup_{\alpha\in \Delta} \|x_\alpha\|_\alpha < \infty\}$ equipped with the norm $\|(x_\alpha)\|:=\sup_{\alpha\in \Delta}\|x_\alpha\|_\alpha.$

Let $f\in L^1(G).$ Then the Fourier transform of $f$ is defined by 
$$\hat{f}(\pi)=\int_G f(x)\pi(x)\, dx, \qquad [\pi]\in \widehat{G}.$$
Let $\mathcal{B}_2(\mathcal{H})$ denote the Hilbert space of all Hilbert-Schmidt operators on a Hilbert space $\mathcal{H}$, with the inner product $\langle T,S\rangle:= tr(TS^\ast).$ Then one can check that $\hat{f}(\pi)\in l^\infty-\bigoplus_{[\pi]\in \widehat{G}}\mathcal{B}_2(\mathcal{H}_\pi)$, where $\mathcal{H}_\pi$ is the Hilbert space corresponding to the representation $\pi$ of $G.$

\begin{theorem}[\cite{Kumar}]  \label{T1.2}
    Suppose that the map $T:L^1(G)\to l^\infty-\bigoplus_{[\pi]\in \widehat{G}} \mathcal{B}_2(\mathcal{H}_\pi)$ is nonzero continuous and $\ast$- preserving and such that:
    \begin{itemize}
        \item [(i)] $T(f\ast g)(\pi)=Tf(\pi)Tg(\pi)$ for all $f,g \in L^1(G);$ and 
        \item [(ii)] $T(R_xf)(\pi)=Tf(\pi)\pi(x)^\ast$ for all $f\in L^1(G)$, $x\in G, [\pi]\in \widehat{G}.$
    \end{itemize}
    Then for each $[\pi]\in E$, we have $T(f)(\pi)= \widehat{f}(\pi)$ for all $f\in L^1(G)$, where 
    $E:=\{[\pi]\in \widehat{G}: Tf(\pi)\neq 0 \text{ for some } f\in L^1(G)\}.$
\end{theorem}
By imposing additional assumption, Krukowski \cite{krukowski2} obtained a much cleaner proof, and the resulting characterization holds on the entire domain. Here we state the Krukowski's result.

\begin{theorem}[\cite{krukowski2}]  \label{T1.3}
    Let $T: L^1(G)\to l^\infty-\bigoplus_{\pi\in \widehat{G}}B(\mathcal{H}_\pi)$ be a linear and bounded operator, which satisfies
    \begin{itemize}
        \item [(i)] $T(f)^\ast= T(f^\ast)$ for every $f\in L^1(G),$
        \item [(ii)] for every $\pi\in \widehat{G}$ the map $f\mapsto T(f)(\pi)$ is irreducible,
        \item [(iv)] for every $f,g \in L^1(G)$ and $\sigma\in \widehat{G}$ we have 
        $$T(f\ast g)(\sigma)=Tf(\sigma)Tg(\sigma),$$
        \item [(iv)] for every $f\in L^1(G),$ $x\in G$ and $\pi\in \widehat{G}$ we have 
        $$T(L_xf)(\pi)= Tf(\pi)\pi(x)^\ast,$$
    \end{itemize}
    then $T$ is the Fourier transform.
\end{theorem}
Note that the above two results are essentially identical. Kumar and Sivananthan do not assume pointwise irreducibility, but their characterization holds only on $E$. Krukowski’s additional assumption yields a simpler proof and a characterization valid on the entire domain. Motivated by the above, we ask whether a similar approach can yield a simpler characterization of the group Fourier transform on the Heisenberg group, valid on the entire domain. We address this in the next section.



The rest of the paper is organized as follows. In Section \ref{section2}, we characterize the group Fourier transform on the Heisenberg group. Sections \ref{section3}–\ref{section6} are devoted to characterizations of $\Hi$-valued integral transforms on locally compact abelian group $G$. In Section \ref{section3}, we characterize the Fourier transform on $L^1(G,\Hi)$, followed by its $L^2(G,\Hi)$ characterization in Section \ref{section4}. In Sections \ref{section5} and \ref{section6}, we consider the cosine and Laplace transforms, respectively.



\section{Group Fourier Transform on the Heisenberg group} \label{section2}
Now we will characterize the group Fourier transform on the Heisenberg group. We aim to extend Krukowski’s result from compact groups to the Heisenberg group, which is locally compact but non-compact, using the same key tool that remains valid for locally compact groups.
Let's begin by recalling the main tool. 
For the following proposition $G$ is any locally compact group.
\begin{propn}[\cite{deitmar2009principles}] \label{P6.1}
    Let $\pi: L^1(G)\longrightarrow \mathcal{B}(V)$ be a nondegenerate $\ast$- representation on a Hilbert space $V$. Then there exists a unique unitary representation $\left(\tilde{\pi}, V\right)$ of $G$ such that $\left\langle\pi(f)v,w\right\rangle= \int_G f(x)\left\langle \tilde{\pi}(x)v, w\right\rangle\, dx$ holds for all $f\in L^1(G)$ and all $v,w\in V.$
\end{propn}

Now we state our result.

\begin{theorem}
    Let $T: L^1(\He)\longrightarrow L^\infty \left(\R^\ast, \mathcal{B}\left(L^2(\R^n)\right), d\mu\right)$ be a bounded linear map which satisfies:
    \begin{itemize}
        \item [(i)] $T(f^\ast)(\lambda)=T(f)(\lambda)^\ast$ for all $\lambda\in \R^\ast, f\in L^1(\He)$,
        \item [(ii)] $T(f\ast g)(\lambda)=Tf(\lambda)Tg(\lambda),$ for all $\lambda\in \R^\ast, f,g \in L^1(\He),$
        \item [(iii)] $T\left(R_{(z,t)}f\right)(\lambda)=Tf(\lambda)\pi_\lambda(z,t)^\ast$
        \item [(iv)] For $\lambda\in\R^\ast$, the map $f\mapsto Tf(\lambda)$ is irreducible.
    \end{itemize}
    Then we have $Tf(\lambda)=\widehat{f}(\lambda)$ for all $\lambda\in \R^\ast$, $f\in L^1(\He).$
\end{theorem}

\begin{proof}
    Let $\lambda\in \R^\ast$, define the map 
    $T_\lambda: L^1(\He)\longrightarrow \mathcal{B}(L^2(\R^n))$ by 
    $T_\lambda(f)= T(f)(\lambda).$ Now using the hypothesis of $T$, we have $T_\lambda$ is a irreducible $\ast$- representation of $L^1(\He).$ Now using the proposition \ref{P6.1}, there exist unique unitary representation $\left(\widetilde{\pi_\lambda}, L^2(\R^n)\right)$ of $\He$ such that $$\left\langle T_\lambda(f)\varphi, \psi\right\rangle =\int_{\He}f(z,t)\left\langle \widetilde{\pi_\lambda}(z,t)\varphi,\psi\right\rangle\, dz\, dt, \quad f\in L^1(\He), ~ \varphi,\psi\in L^2(\R^n).$$
    Now let $(w,s)\in \He$, $\pi_\lambda$ be the Schr\"odinger representation of $\He$ and compute,
    \begin{align*}
        \int f(z,t)\left\langle \widetilde{\pi}_\lambda(z,t)\pi_\lambda(w,s)^\ast \varphi, \psi\right\rangle \, dz\, dt &= \int f(z,t)\left\langle \widetilde{\pi_\lambda}(z,t)\pi_\lambda(w,s)^\ast \varphi, \psi\right\rangle\, dz\, dt\\
        &= \left\langle T_\lambda(f)\pi_\lambda(w,s)^\ast \varphi, \psi\right\rangle\\
        &= \left\langle T\left( R_{(w,s)}f\right)(\lambda)\varphi, \psi\right\rangle\\
        &=\int R_{(w,s)}f(z,t)\left\langle \widetilde{\pi}_\lambda(z,t)\varphi,\psi\right\rangle\, dz\, dt\\
        &= \int f(z,t)\left\langle \widetilde{\pi}_\lambda\left(\left(z,t\right)\left(w,s\right)^{-1}\right)\varphi, \psi\right\rangle\, dz\, dt\\
        &= \int f(z,t)\left\langle \widetilde{\pi}_\lambda(z,t)\widetilde{\pi}_\lambda(w,s)^\ast \varphi, \psi\right\rangle\, dz\, dt.
    \end{align*}
    This implies 
    $$\int f(z,t)\left\langle \widetilde{\pi}_\lambda(z,t)\left(\pi_\lambda(w,s)^\ast-\widetilde{\pi}_\lambda(w,s)^\ast\right)\varphi, \psi\right\rangle\, dz\, dt=0,$$
    which is true for all $f\in L^1(\He),~(z,t)\in \He$ and $\varphi,\psi\in L^2(\R^n)$. Now if we put $$f(z,t)=\overline{\left\langle \widetilde{\pi}_\lambda(z,t)\left(\pi_\lambda(w,s)^\ast-\widetilde{\pi}_\lambda(w,s)^\ast\right)\varphi,\psi\right\rangle },$$ we obtain $$\left\langle \widetilde{\pi}_\lambda(z,t)\left(\pi_\lambda(w,s)^\ast-\widetilde{\pi}_\lambda(w,s)^\ast\right)\varphi,\psi\right\rangle=0,$$ for all $(z,t), (w,s)\in \He$ and $\varphi,\psi\in L^2(\R^n).$ Finally taking $(z,t)$ to be the identity element of $\He$ and $\psi=\left(\pi_\lambda(w,s)^\ast-\widetilde{\pi}_\lambda(w,s)^\ast\right)\varphi$, we get $\left(\pi_\lambda(w,s)^\ast-\widetilde{\pi}_\lambda(w,s)^\ast\right)\varphi=0$ for all $(w,s)\in \He$ and $\varphi\in L^2(\R^n)$. Which implies $\widetilde{\pi}_\lambda=\pi_\lambda.$
    Therefore we have 
    $$\left\langle T(f)(\lambda)\varphi, \psi\right\rangle= \int f(z,t)\left\langle \pi_\lambda(z,t)\varphi, \psi\right\rangle\, dz\, dt.$$
    Hence $T$ is nothing but the group Fourier transform, i.e. $Tf(\lambda)=\widehat{f}(\lambda) ~$ for all $f\in L^1(\He)$ and $\lambda\in \R^\ast.$
\end{proof}

\begin{remark}
    We note that in this theorem we have used extra assumption of irreducibility as compared to the theorem \ref{T1.1} by Lavanya and Thangavelu \cite{Lavanya_Heisen}, but we have obtain the characterization on the entire domain. The main crux of the proof is the bijection between the unitary representations of $\Hi$ and irreducible, unitary $*$-representations of $L^1(\Hi)$, which holds for any locally compact group $\Hi$. Arguably, our proof simplifies earlier works of Lavanya and Thangavelu \cite{Lavanya_Heisen}.
\end{remark}

\section{ Fourier Transform on a locally compact abelian group} \label{section3}

\noindent
In this section, we will see a characterization of the 
Fourier transform on a locally compact abelian group for a Hilbert space valued function. Biswas and Radha \cite{Hasan1} characterized the Hilbert space valued Fourier transform on the Euclidean spaces. Our result will generalize their work to the locally compact abelian group. Also our result is much stronger, because we need less hypothesis to characterize the Fourier transform. Here, we adapt the technique of Krukowski to characterize the Fourier transform. Now we discuss some preliminaries for the Fourier transform of locally compact abelian group.

Let $G$ be a locally compact Abelian group with Haar measure $\mu$ and $\widehat{G}$ be the dual group. Operations on $G$ will be denoted additively. The Fourier transform $\F: L^1(G)\to C_0(\widehat{G})$ is defined by
$$\F(f)(\chi)= \int_G f(\chi)\overline{\chi(x)}\, d\mu(x).$$
\noindent
The convolution $\ast: L^1(G)\times L^1(G)\to L^1(G)$ can be defined as 
$$f\ast g(x)=\int_G f(u)g(x-u)\, d\mu(x).$$
Then it is easy to check that $\F(f\ast g)= \F(f)\F(g)$ for all $f,g\in L^1(G).$
Now we state the Result of Krukowski \cite{krukowski1} for characterization of Fourier transform for Locally compact abelian groups.
\begin{theorem}
    Let $T: L^1(G)\to C^b(\widehat{G})$ be a linear and bounded transform. If it satisfies $T(f\ast g)= T(f)T(g)$ for all $f,g \in L^1(G)$. Then there exist an open set $E\subset \widehat{G}$ and a continuous function $\Theta: E\to \widehat{G}$ such that 
    $$Tf(\phi)=\begin{cases}
        \F f\circ \Theta(\phi) & \text{ if } \phi \in E\\
        0 & \text{ if } \phi \in \widehat{G}\setminus E.
    \end{cases}$$
    Furthermore, if we assume $T$ satisfies the additional property $T(L_xf)(\phi)= \overline{\phi(x)}T(f)(\phi)$ for all $f\in L^1(G), \phi \in \widehat{G}, x\in G$, then $\Theta(\phi)= \Bar{\phi}.$
\end{theorem}

Note that for $\phi \in E,$ we choose $g_{\ast}\in L^1(G)$ such that $1=T(g_{\ast})(\phi).$
Now, we can reinterpret his results as follows.
For $\phi \in E,$ we have
\begin{align*}
  &T(f)(\phi)=T(f)(\phi) T(g_\ast)(\phi)=T(f\ast g_\ast)(\phi) =T_\phi\left(\int_G f(y)L_y g_\ast \, dy\right) \\
  &=\int_G f(y)T_\phi\left(L_y g_\ast\right)\, dy= \int_G f(y)\overline{\phi(y)}T(g_\ast)(\phi)\, dy=\int f(y)\overline{\phi(y)}\, dy=\F(f)(\phi).
\end{align*}
This observation helps us to extend this result in the context of a Hilbert space valued function, see \ref{T3.9} below.
In the paper \cite{Hasan1}, Biswas and Radha have characterized the Hilbert space valued Fourier transform on Euclidean space. Throughout this paper, we denote by $(\Hi, \langle {\cdot, \cdot } \rangle)$, a separable Hilbert space with an orthonormal basis $\{e_k: k\in\N\}$.   Then the $\ast$- convolution of two $\Hi$- valued function is defined in \cite{Hasan1} as follows.
\begin{equation} \label{eq: vector conv def}
    (f\ast g) (x)=\int_{\R}\left\langle f(y), g(x-y)\right\rangle dy,~~f,g\in L^1(\R,\Hi).
\end{equation}
It is easy to see that for $f,g\in L^1(\R,\Hi)$, $f\ast g\in L^1(\R)$. So, the authors in \cite{Hasan1} defined the $\odot$- convolution of a scalar valued function and a vector valued function as follows.
\begin{equation}\label{eq: vector scalar conv}
    (f\odot g)(x)=\int_\R f(x-y)g(y)dy,~~f\in L^1(\R), g\in L^1(\R,\Hi).
\end{equation}
Then it is easy to observe that $(f\ast g)\odot h\in L^1(\R, \Hi)$ whenever $f,g,h\in L^1(\R,\Hi)$.
Now the Fourier transform of $f\in L^1(\R,\Hi)$ is defined by 
\begin{equation*}
    \widehat{f}(\xi)=\int_\R e^{-2\pi ix\xi}f(x)dx,~~\xi\in R,
\end{equation*}
where the right hand side is a Bochner integral. For a detailed study of vector-valued integrations, we refer to \cite{Diestel}.
We now state the characterization result of Biswas and Radha \cite{Hasan1} for vector-valued functions.

\begin{theorem} \label{L1 char}
     Let $(\Hi,\langle \cdot, \cdot \rangle)$ be a separable Hilbert space with countable orthonormal basis $\{e_1,e_2, \cdot \cdot \cdot\}$ and $T: L^1(\R,\Hi) \to C_0(\R,\Hi)$ be a bounded linear map satisfying:
        \begin{enumerate}
            \item $ T((f*g)\odot h)(\xi)= \langle Tf(\xi),Tg(-\xi)\rangle Th(\xi)$, for $f,g,h \in L^1(\R,\Hi), \xi \in \R$.
            \item  $T(\tau_s f)(\xi)= e^{-2\pi i s \xi}Tf(\xi), f \in L^1(\R,\Hi), s,\xi \in \R$.
            \item If $\langle f(x),e_k \rangle =0 $ for a.e. $x \in \R$ then $\langle Tf(\xi), e_k \rangle =0 $ for all $\xi \in \R$.
            \item $T(\Tilde{f})(-\xi)=Tf(\xi), f \in L^1(\R,\Hi),\xi \in \R$, where $\Tilde{f}(x)=f(-x)$.
        \end{enumerate}
        Then there exists $E \subset \R$ such that 
        \begin{center}
            $\langle Tf(\xi),Tg(\xi)\rangle = \chi_E(x)\langle \widehat{f}(\xi), \widehat{g}(\xi)\rangle$  $f,g \in L^1(\R,\Hi), \xi \in \R$
        \end{center}
\end{theorem}


Notice that, as a corollary to Theorem \ref{L1 char}, they recover the Fourier transform.
\begin{corollary}[\cite{Hasan1}]\label{cor: L^1 cor}
If in addition to the hypothesis of Theorem \ref{L1 char}, $T$ satisfies $T(g_k)=g_k$ where $g_k(x)=e^{-\pi x^2}e_k$ for all $k\in \N$, then there exists $E\subset \R$ such that $Tf(\xi)=\chi_E (\xi) \widehat{f}(\xi)$,\quad $f \in L^1(\R, H),~~\xi \in \R$.
\end{corollary}
Now we wish to get analogous results for locally compact Abelian groups. For that at first we need to defined corresponding Fourier transform and convolutions. Let $G$ be a locally compact abelian group with the Haar measure $\mu$. For simplicity we write $dx$ instead of $d\mu(x).$

\begin{definition}
    Let $f\in L^1(G,\Hi)$, then the Fourier transform of $f$ is defined by
    $$\F f(\chi)= \int_G \overline{\chi(x)} f(x) \, dx, \qquad ~ \text{ for } \chi\in \hat{G},$$
    where the right hand side is a Bochner integral.
\end{definition}
Then it is easy to see that $\F f\in C_0(\widehat{G},\Hi).$
The convolution of $\Hi$-valued function can be defined as follows:
\begin{equation}
    f\ast g(x)= \int_G \left\langle f(y), g(x-y)\right\rangle \, dy, ~~ f,g \in L^1(G,\Hi).
\end{equation}
Observe that for $f,g\in L^1(G,\Hi)$, $f\ast g\in L^1(G,\C)$. After that as in \cite{Hasan1}, the convolution of a scalar valued function and a vector valued function can be defined as follows.
\begin{equation} \label{E3.4}
    f\odot g(x)=\int_G f(y) g(x-y)\, dy,~~ f\in L^1(G,\C), ~g\in L^1(G,\Hi). 
\end{equation}
    Then we can observe that $f\odot g \in L^1(G,\Hi)$. Therefore we have $(f\ast g)\odot h \in L^1(G, \Hi)$ whenever $f,g,h \in L^1(G,\Hi).$
 
\noindent
Notation: We use both $\hat{f}$ and $\F(f)$ to denote the Fourier transform. For the translation operator, we use $L_x$ and $\tau_x$ interchangeably, with
$
\tau_x f(y)=f(y\cdot x^{-1}),
$
where $\cdot$ denotes the group operation.


Next, we have the following convolution theorem in this context.

\begin{theorem}[Fourier Convolution Property]
    Let $f,g,h \in L^1(G,\Hi)$. Then we have $$\F \left(\left(f\ast g\right)\odot h\right)\left(\chi\right)= \left\langle \F f\left(\chi\right), \F g\left(\bar{\chi}\right)\right\rangle \F h\left(\chi\right)$$
\end{theorem}

\begin{proof}
    Let $f,g\in L^1(G,\Hi)$.
    \begin{align*}
        &\F\left(f\ast g\right)(\chi)
        = \int\left(f\ast g\right)\overline{\chi(x)}\, dx= \int \int \left\langle f(y), g(x-y)\right\rangle\overline{\chi(x)}\, dy\, dx\\
        &= \int \left\langle f(y), \int g(x-y)\chi(x)\, dx\right\rangle\, dy\\
        &= \int \left\langle f(y),\int g(x)\chi(x+y)\, dx\right\rangle\, dy\\
        &=\int \left\langle f(y)\overline{\chi(y)}, \F g(\overline{\chi})\right\rangle\, dy\\
        &= \left\langle \F f(\chi), \F g(\overline{\chi})\right\rangle.
    \end{align*}
    Now if $\phi\in L^1(G,\C)$ and $h\in L^1(G,\Hi)$, then
    \begin{align*}
        \F(\phi\odot h)(\chi)
        &= \int \phi\odot h(x)\overline{\chi(x)}\, dx\\
        &=\int \int \phi(y)h(x-y)\overline{\chi(x)}\, dx\, dy\\
        &=\int \phi(y)\left(\int h(x)\overline{\chi(x+y)}\, dx\right)\, dy\\
        &= \left(\int \phi(y)\overline{\chi(y)}\, dy\right)
        \left(\int h(x)\overline{\chi(x)}\, dx\right)\\
        &=\F \phi(\chi)\F h(\chi).
    \end{align*}
    Since $f\ast g\in L^1(G,\C)$ and $h\in L^1(G,\Hi)$, combining the above two, we can conclude that
    $$\left(\left(f\ast g\right)\odot h\right)\left(\chi\right)= \left\langle \F f\left(\chi\right), \F g\left(\bar{\chi}\right)\right\rangle \F h\left(\chi\right).$$
\end{proof}
  Now for $s\in G$ the translation operator $L_s$ can be defined by $L_s f(g)= f(gs^{-1})$ for all $g\in G.$ Now we will see how Fourier transform will act on translation operator. A simple change of variable gives the following proposition.
  \begin{propn}
      Let $f \in L^1(G,\Hi)$ and $s\in G$. Then we have $\F (L_s f)(\phi)= \overline{\phi(s)}\F f(\phi)$, for all $\phi \in \widehat{G}.$
  \end{propn}
Now we need to state a lemma which will be needed for our proof.
\begin{lemma}\cite{Inf_Dim_Ana} \label{L3.7}
    Let $f: \Omega \to X$ be Bochner integrable and $Y$ be a Banach space. If $T: X \to Y$ is a bounded operator, then the function $Tf: \Omega \to Y$ defined by $Tf(\omega)= T(f(\omega))$ is also Bochner integrable and $$\int_{\Omega}Tf d\mu= T\left( \int_{\Omega}f d\mu\right).$$
\end{lemma}

\begin{lemma} \label{L3.8}
    For $f,g,h\in L^1(G,\Hi)$ and $a,b\in G$ the following equality holds
    $$\left(L_af\ast g\right)\odot L_b h= \left(f\ast L_{a+b}g\right)\odot h$$
\end{lemma}

\begin{proof}
A simple change of variables yields.
    \begin{align*}
       \left(L_af\ast g\right)\odot L_b h(x) &= \int_G \int_G \left\langle L_a f(z), g(y-z)\right\rangle L_b h(x-y)\, dy\, dz\\
       &= \int_G\int_G \left\langle f(z-a), g(y-z)\right\rangle h(x-y-b)\, dy\, dz\\
       &=\int_G \int_G \left\langle f(z), g(y-z-a-b)\right\rangle h(x-y)\, dy\, dz\\
       &= \left(f\ast L_{a+b} g\right)\odot h(x).
    \end{align*}
\end{proof}

\begin{remark}
    For our proof we also need a orthogonality condition. Indeed for the Fourier transform one can easily check that, if $\langle f(x),e_k\rangle=0$ for a.e. $x\in G$, then we have $\langle \widehat{f}(\phi), e_k\rangle =0$ for all $\phi\in \widehat{G}.$ One can also check that the Fourier transform also satisfies $\F\left(\widetilde{f}\right)\left(\overline{\phi}\right)=\F(f)(\phi).$
\end{remark}
Now we state our result.
\begin{theorem} \label{T3.7}
    Let $G$ be a locally compact abelian group. $T:L^1(G,\Hi)\to C_0(\widehat{G},\Hi)$ be a bounded linear map satisfying 
    \begin{itemize}
        \item [(i)] $T\left(\left(f\ast g\right)\odot h\right)(\phi)=\left\langle Tf(\phi),Tg\left(\overline{\phi}\right)\right\rangle Th(\phi)$ \quad
        for all $f,g,h \in L^1(G,\Hi), \phi\in \widehat{G}.$
        \item [(ii)] If for some $k$, $\langle f(x), e_k\rangle=0$ for a.e. $x\in G$, then $\langle Tf(\phi), e_k\rangle=0$\quad for all $\phi \in \widehat{G}.$
        \item [(iii)] $T\left( \widetilde{f}\right)\left(\overline{\phi}\right)=T(f)(\phi)$ \quad for all $h\in L^1(G,\Hi)$ and $\phi\in \widehat{G}.$
    \end{itemize}
    Then there exist an open set $E\subset \widehat{G}$ and a continuous map $\beta: E\to \widehat{G}$ such that 
    $$\left\langle Tf\left(\phi\right), Tg\left(\overline{\phi}\right)\right\rangle=
    \begin{cases}
        \left\langle \widehat{f}\left(\beta{(\phi)}\right), \widehat{g}\left(\overline{\beta\left({\phi}\right)}\right)\right\rangle & \text{ if } \phi \in E\\
        0  &\text{ if } \phi \notin E.
    \end{cases}$$
    for all $f,g,h\in L^1(G,\Hi), ~ \phi\in \widehat{G}.$
\end{theorem}

\begin{proof}
    Let $E:=\left\{ \phi\in \widehat{G}: T(f)(\phi)\neq 0 \text{ for some } f\in L^1(G,\Hi)  \right\}$. Then it is easy to see that $E$ is an open set. If $\phi\notin E$, there is nothing to show. So let $\phi\in E$. Then there exist $f\in L^1(G,\Hi)$ such that $T(f)(\phi)\neq 0$. Now define the map $T_\phi: L^1(G,\Hi)\to \Hi$ by $T_\phi(g)= T(g)(\phi)$. It is a bounded linear map and also nonzero. Then by the orthogonality condition $(ii)$, we can get a $h\in L^1(G,\Hi)$ such that $T_\phi(h)=e_k.$ Now define the function $\chi_\phi: G\to \C$ by
    $$\chi_\phi(x):= \left\langle T_\phi\left(L_x h\right), e_k\right\rangle.$$
    Note that the choice of the function $h$ is not unique, so there might be different $\chi_\phi$ corresponding to $\phi$.

    We will now show that $\chi_\phi$ is a homomorphism with range in $S^1.$
    Using Lemma \ref{L3.8}, we obtain
    \begin{align*}
        &\chi_\phi(x)\chi_\phi(y) = \left\langle T\phi(L_x h), e_k\right\rangle \left\langle  T\phi(L_y h), e_k\right\rangle= \left\langle  \left\langle T_\phi(L_x h), e_k\right\rangle T_\phi(L_y h), e_k\right\rangle\\
        &= \left\langle  \left\langle T_\phi(L_x h), T_\phi(h)\right\rangle T_\phi(L_y h), e_k\right\rangle=\left\langle  \left\langle T_\phi(L_x h), T_{\overline{\phi}}(\tilde{h})\right\rangle T_\phi(L_y h), e_k\right\rangle\\
        &= \left\langle T_\phi\left(\left(L_x h\ast \tilde{h}\right)\odot L_y h\right), e_k\right\rangle= \left\langle T_\phi\left(\left( h\ast L_{x+y}\tilde{h}\right)\odot h\right), e_k\right\rangle\\
        &=\left\langle \left\langle T(h)(\phi), T\left(L_{x+y}\tilde{h}\right)(\overline{\phi})\right\rangle T(h)(\phi), e_k\right\rangle\\
        &= \left\langle \left\langle e_k, T\left(L_{-x-y}h\right)(\phi)\right\rangle e_k,e_k\right\rangle\\
        &= \left\langle e_k, T\left(L_{-x-y}h\right)(\phi)\right\rangle = \overline{\chi_{\phi}(-x-y)}.
    \end{align*}
    Now we see that $\chi_\phi(0)=\left\langle T_\phi\left(L_0 h\right),e_k\right\rangle=1.$ Now if we put $y=0$ in the above identity, we obtain $\chi_\phi(x)\chi_\phi(0)=\overline{\chi_\phi(-x)}$, i.e. $\overline{\chi_\phi(-x)}= \chi_\phi(x).$ Hence we get $\chi_\phi(x+y)=\chi_\phi(x)\chi_\phi(y).$ Now we will show that $|\chi_\phi(x)|=1.$ For that compute
    \begin{align*}
       & \chi_\phi(x)\overline{\chi_\phi(x)} = \left\langle T_\phi(L_x h), e_k\right\rangle \left\langle e_k, T_\phi(L_x h)\right\rangle= \left\langle e_k, \left\langle e_k, T_\phi(L_x h)\right\rangle T_\phi(L_x h)\right\rangle\\
        &= \left\langle e_k, \left\langle T_\phi(h),T_{\overline{\phi}}(L_{-x}\tilde{h})\right\rangle T_\phi(L_x h)\right\rangle= \left\langle e_k, T_\phi\left( \left( h\ast L_{-x}\tilde{h}\right)\odot L_xh\right)\right\rangle\\
        &= \left\langle e_k, T_\phi\left(\left( h\ast L_{x}L_{-x}\tilde{h}\right)\odot h\right)\right\rangle= \left\langle e_k, T_\phi\left(\left( h\ast\tilde{h}\right)\odot h\right)\right\rangle\\ 
        &= \left\langle e_k, \left\langle T_\phi(h), T_{\overline{\phi}}(\tilde{h})\right\rangle T_\phi(h)\right\rangle\\
        &= \left\langle e_k, \left\langle e_k, e_k \right\rangle e_k\right\rangle \\
        &= 1.
    \end{align*}
    Now using the continuity of $T_\phi$ and by \cite[lemma 1.4.2]{deitmar2009principles}, we can get $\chi_\phi$ is continuous. Therefore $\chi_\phi\in \widehat{G}.$ Now we set $\beta(\phi):= \overline{\chi_\phi}$ and using lemma \ref{L3.7} (see also, \cite[Proposition 7]{Vec_Meas}) we have
    \begin{align*}
       & \left\langle T(f)(\phi), T(g)(\overline{\phi})\right\rangle
        =  \left\langle T(f)(\phi), T(g)(\overline{\phi})\right\rangle \left\langle T(h)(\phi), e_k\right\rangle=\left\langle  \left\langle T(f)(\phi), T(g)(\overline{\phi})\right\rangle T(h)(\phi), e_k\right\rangle\\
        &= \left\langle T\left(\left(f\ast g\right)\odot h\right)(\phi), e_k\right\rangle= \left\langle T_\phi\left(\int_G\int_G \left\langle f(z), g(y-z)\right\rangle L_y h\, dy\, dz\right), e_k\right\rangle\\
        &= \left\langle \int_G\int_G \left\langle f(z), g(y-z)\right\rangle T_\phi\left(L_y h\right)\, dy\, dz, e_k\right\rangle=  \int_G\int_G \left\langle f(z), g(y-z)\right\rangle \left\langle T_\phi\left(L_y h\right), e_k\right\rangle\, dy\, dz\\
        &=  \int_G\int_G \left\langle f(z), g(y-z)\right\rangle \chi_\phi(y), dy\, dz=  \int_G\int_G \left\langle f(z), g(y)\right\rangle \chi_\phi(y+z), dy\, dz\\
        &=  \int_G\int_G \left\langle f(z), g(y)\right\rangle \chi_\phi(y)\chi_\phi(z), dy\, dz=  \int_G\int_G \left\langle f(z)\chi_\phi(z), g(y)\overline{\chi_\phi(y)}\right\rangle , dy\, dz\\
        &=  \left\langle \int_G f(z)\chi_\phi(z)\, dz, \int_G g(y)\overline{\chi_\phi(y)}\, dy\right\rangle = \left\langle \widehat{f}\left(\overline{\chi_\phi}\right), \widehat{g}\left(\chi_\phi\right)\right\rangle= \left\langle \widehat{f}\left(\beta\left(\phi\right)\right), \widehat{g}\left(\overline{\beta\left(\phi\right)}\right)\right\rangle.
    \end{align*}
It remains only to show that $\beta$ is continuous to complete the proof.
 using \cite[lemma 1.4.2]{deitmar2009principles} and continuity of $T$, we have $x\mapsto T\left(L_x h\right)$ is continuous. Therefore the family $\left(T\left(L_x h\right)\right)_{x\in K}$ is compact in $C^b\left(\widehat{G}, \Hi\right)$ for every compact $K\subset G.$ Now by Arzel\'a-Ascoli theorem, the family is equicontinuous at every $\phi'\in E$, i.e., for every $\epsilon>0$ there exists an open neighbourhood $U'$ of $\phi'$ such that
$$\|T\left(L_x h\right)(\phi)-T\left(L_x h\right)(\phi')\|_H<\epsilon$$ 
for every $\phi\in U'$ and $x\in K.$ Now if we calculate 
\begin{align*}
    |\overline{\chi_\phi(x)}-\overline{\chi_{\phi'}(x)}|
    &= |\left\langle e_k, T(L_x h)(\phi)\right\rangle-\left\langle e_k, T(L_x h)(\phi')\right\rangle|\\
    &= |\left\langle e_k, T(L_x h)(\phi)- T(L_x h)(\phi')\right\rangle|\\
    &\leq \|e_k\|_H\|T(L_x h)(\phi)- T(L_x h)(\phi')\|_H.
\end{align*}
    Therefore we have 
    $$\sup_{x\in K}|\overline{\chi_\phi(x)}-\overline{\chi_{\phi'}(x)}|\leq \epsilon$$
    for every $\phi\in U'$. This gives that if $\phi$ converges to $\phi'$, $\overline{\chi_\phi}$ converges to $\overline{\chi_{\phi'}}$. Hence $\beta$ is continuous. This completes the proof of the theorem.
\end{proof}
Note that under the additional assumption: $T(L_s h)(\phi)= \overline{\phi(s)}T(h)(\phi),$ it is easy to see that $\beta(\phi)=\phi$. In the next theorem, we shall give a simpler proof of characterization of the Fourier transform on $L^1(G,\Hi).$ The proof of the above theorem shows that the translation–modulation condition itself is sufficiently stronger. Note that Biswas and Radha \cite{Hasan1} in Theorem \ref{L1 char} assumed this condition together with two additional hypothesis. However, the following theorem shows that the translation–modulation condition alone suffices. 
Thus, our result provides a stronger characterization.


\begin{theorem} \label{T3.9}
    Let $G$ be a locally compact abelian group. $T: L^1(G, \Hi)\to C_0(\widehat{G}, \Hi)$ be a bounded linear map satisfying the following:
    \begin{itemize}
        \item [(i)] $T\left(\left(f\ast g\right)\odot h\right)(\phi)= \left\langle Tf(\phi), Tg(\Bar{\phi})\right\rangle Th(\phi)$
        \item [(ii)] $T(L_s h)(\phi)= \overline{\phi(s)}T(h)(\phi).$
    \end{itemize}
    Then there exist $E\subset \widehat{G}$ such that 
    
    $$\left\langle Tf\left(\phi\right), Tg\left(\overline{\phi}\right)\right\rangle=
    \begin{cases}
        \left\langle \widehat{f}\left(\phi\right), \widehat{g}\left(\overline{\phi}\right)\right\rangle & \text{ if } \phi \in E\\
        0  &\text{ if } \phi \notin E.
    \end{cases}$$
    for all $f,g,h\in L^1(G,\Hi), ~ \phi\in \widehat{G}.$
\end{theorem}
\begin{corollary}
If in addition to the hypothesis of Theorem \ref{T3.9}, $T$ satisfies $T(g_k)=\widehat{g}(\phi)e_k$, for all $k\in \N$ where $g_k=ge_k$ for some $g\in L^1(G)$, with $\hat{g}(\phi)\neq 0$, for all $\phi\in \widehat{G}$.
    Then we have $Tf(\phi)=\widehat{f}(\phi)$ for all $\phi \in E.$
\end{corollary}

\begin{proof}[Proof of Theorem \ref{T3.9}.]
    Let $E= \left\{\phi \in \widehat{G}: T(f)(\phi)\neq 0, \text{ for some }f \in L^1(G,\Hi)\right\}$. If $\phi \in E^c$, then we are done. If $\phi \in E$, define $T_\phi: L^1(G,\Hi)\to \Hi$ by $T_\phi(f):= T(f)(\phi)$. Then $T_\phi$ is linear. Now $\phi \in E$, implies there exists $h\in L^1(G,\Hi)$ such that $T(h)(\phi)\neq 0$, which means for some $k\in \N$, $\langle T(h)(\phi), e_k\rangle \neq 0$. WLOG we can take $\langle T(h)(\phi), e_k \rangle =1$. Now we compute the following:
    \begin{align*}
       & \left\langle Tf(\phi), Tg(\bar{\phi})\right\rangle
        = \left\langle Tf(\phi), Tg(\bar{\phi})\right\rangle \left\langle T(h)(\phi), e_k \right\rangle = \left\langle \left\langle Tf(\phi), Tg(\bar{\phi}) \right\rangle T(h)(\phi), e_k \right\rangle\\
        &= \left\langle T\left(\left( f\ast g\right)\odot h\right)(\phi), e_k\right\rangle= \left\langle T_\phi\left(\left( f\ast g\right)\odot h\right), e_k\right\rangle\\
        &= \left\langle T_\phi \left( \int \int \left\langle f(z), g(y-z)\right\rangle L_y(h) \, dy\, dz\right), e_k\right\rangle \\
        &= \left\langle \int \int \left\langle f(z), g(y-z)\right\rangle  T_\phi\left(L_y(h)\right) \, dy\, dz, e_k\right\rangle \\
        &= \left\langle \int \int \left\langle f(z), g(y-z)\right\rangle  \overline{\phi(y)}T_\phi(h) \, dy\, dz, e_k\right\rangle \\
        &= \left\langle \int \int \left\langle f(z), g(z)\right\rangle  \overline{\phi(y)}\hspace{2pt}\overline{\phi(z)}T_\phi(h) \, dy\, dz, e_k\right\rangle \\
        &= \left\langle \left\langle \F f(\phi), \F g(\Bar{\phi})\right\rangle T_\phi(h), e_k\right\rangle\\
        &= \left\langle \F f(\phi), \F g(\Bar{\phi})\right\rangle \left\langle T_\phi(h), e_k\right\rangle\\
        &= \left\langle \F f(\phi), \F g(\Bar{\phi})\right\rangle.
        \end{align*}
        So, we have obtain $\left\langle Tf(\phi), Tg(\Bar{\phi})\right\rangle= \chi_E(\phi)\left\langle \F f(\phi), \F g(\Bar{\phi})\right\rangle.$
\end{proof}

\begin{remark}
    Adapting Krukowski's technique \cite{krukowski1}, our proof gives an easy and compact proof comparing to Hasan-Radha's proof \cite{Hasan1} for characterization on $L^1(\R,\Hi)$. The most important thing is that we don't assume the orthogonality condition and the condition $T(\widetilde{f})=\widetilde{Tf}$, as in Theorem \ref{L1 char}.
\end{remark}

\section{Hilbert Space Valued $L^2$ Characterization} \label{section4}

In this section, we characterize the Fourier transform in $L^2(G,\Hi)$, where $G$ is a locally compact abelian group and $\Hi$ is a Hilbert $A$- module for a $C^\ast$- algebra $A$. For a study of Hilbert $C^\ast$-modules, we refer \cite{Lance_Book}. For $L^2(G,\Hi)$, we note that our result extends to general Hilbert $A$-modules.
We now discuss the main ingredients required for the proof of our result.


\begin{theorem}\cite{Larsen_Mult_Book}  \label{T4.1}
   Let $G$ be a locally compact abelian group and suppose $T: L^2(G) \to L^2(G) $ is a linear transformation. Then the following are equivalent:
   \begin{itemize}
       \item [(i)] $T\in M\left(L^2(G)\right)$, i.e., $T$ is translation-invariant.
       \item [(ii)] There exists a unique $\varphi\in L^\infty (\widehat{G})$ such that $\left(Tf\right)^{\widehat{}}=\varphi \widehat{f}$ for each $f\in L^2(G).$
   \end{itemize}
    Moreover the correspondence between $T$ and $\varphi$ defines an isometric algebra isomorphism from $M\left(L^2(G)\right)$ onto $L^\infty(\widehat{G}).$
\end{theorem}

\noindent
Now we will prove the vector valued analogue of the above result. 
\begin{theorem} \label{T4.2}
    Let $T: L^2(G,\Hi)\to L^2(G,\Hi)$ be a bounded linear operator of the form $T= S\otimes I$, where $S\in B(L^2(G))$ such that $T$ commutes with translation, then there exists $m\in L^\infty(G)$, such that $(Tf)~\widehat{}~(\chi)= m(\chi)\widehat{f}(\chi)$, $f\in L^2(G,\Hi)$, $\chi \in \widehat{G}.$
\end{theorem}

\begin{proof}
    It is enough to prove it for basic tensors. Let $f= F\otimes u$ be a basic tensor, $f\in L^2(G)$, $u\in \Hi.$ Then $T(L_t(F\otimes u))= \tau_t((S\otimes I)(F\otimes u))= \tau_t(SF \otimes u)= \tau_t SF\otimes u.$ On the other hand, $T(\tau_t(F\otimes u))= (S\otimes I)(\tau_t F\otimes u))= (S\tau_tF)\otimes u.$ Hence we have, $S\tau_t F\otimes u= \tau_t SF\otimes u$, implies, $(S\tau_t F- \tau_t SF)\otimes u=0$. Again which implies $(S\tau_t F- \tau_t SF)\otimes \langle u, v\rangle =0$ for all $v\in \Hi.$ Now if we choose $v$ such that $\langle u,v \rangle \neq 0$, we obtain $(S\tau_t F = \tau_t SF).$ Hence $S: L^2(G)\to L^2(G)$ commutes with translation. Then by theorem \ref{T4.1}, there exists $m\in L^\infty(\widehat{G})$, such that $(Sf)~\widehat{}~(\chi)= m(\chi)\widehat{f}(\chi)$. Thus we have 
    $$(Tf)~\widehat{}~(\chi)= (SF\otimes u)~\widehat{}~(\chi)= (SF)~\widehat{}~(\chi)\otimes u= m(\chi) \widehat{F}(\chi)\otimes u = m(\chi) \widehat{f}(\chi),$$ which proves our assertion.
\end{proof}
Now we state our result.
\begin{theorem}
Let $\left(\Hi, \left\langle \cdot , \cdot \right\rangle\right)$ be a Hilbert A-module with countable orthonormal basis $\{e_1, e_2, \ldots\}$ and $T : L^2\left(G, \Hi\right) \to L^2(\widehat{G}, \Hi)$ be a bounded, linear and surjective map satisfying:
\begin{itemize}
    \item [(i)] $T ((f \ast g)\odot h)(\chi ) = \left\langle T f (\chi ), T g(\overline{\chi} ) \right\rangle T h(\chi),$\quad for $f, g, h \in L^2(G, \Hi)$, $(f \ast g)\odot h \in L^2(G, \Hi)$, $\chi \in \widehat{G}.$
    \item [(ii)] $T (\tau_s f )(\chi ) = \overline{\phi(s)} T f (\chi ),$\quad $f \in L^2(G, \Hi)$, $s \in G$ and $\chi\in \widehat{G}.$
    \item [(iii)] $T$ is of the form $S \otimes I$ for some $S \in \mathcal{B}\left(L^2(G),L^2(\widehat{G})\right)$, $I$ is the identity operator on $\Hi$.
    \item [(iv)] $T (\tilde{f })(\overline{\phi} ) = T f (\phi ),$ \quad $f \in L^2(G, \Hi)$, $\phi \in \widehat{G}$, where $\tilde{f }(x) = f (-x).$
\end{itemize}
Then $$\langle T f (\xi ), T g(\xi ) \rangle = \langle \widehat{f }(\xi ), \widehat{g}(\xi ) \rangle, \quad f, g \in L^2(G, \Hi), \xi \in G.$$
\end{theorem}

\begin{proof}
    We have $T: L^2(G,\Hi)\to L^2(\widehat{G},\Hi)$ be the map satisfying those properties. Now using it define a map $U: L^2(G,\Hi)\to L^2(G,\Hi)$ by $Uf= \F^{-1}(Tf)$. Then clearly $U$ is linear, bounded and surjective. Now let $f_1,f_2, f_3\in L^2(G,\Hi)$ such that $(f_1\ast f_2)\odot f_3\in L^2(G,\Hi).$ Let $Uf_1=g_1, Uf_2=g_2, Uf_3=g_3$ and $U(f_1\ast f_2)\odot f_3= g$. Then $Tf_1= \widehat{g_1}, Tf_2=\widehat{g_2}, Tf_3= \widehat{g_3}$ and $T(f_1\ast f_2)\odot f_3= \widehat{g}.$ Now for $\chi\in \widehat{G}$, we have
    \begin{align*}
        \widehat{g}(\chi) &= T((f_1\ast f_2)\odot f_3)(\chi)\\
        &= \left\langle Tf_1(\chi), Tf_2(\bar{\chi})\right\rangle Tf_3(\chi)\\
        &= \left\langle \widehat{g_1}(\chi), \widehat{g_2}(\bar{\chi})\right\rangle \widehat{g_3}(\chi)\\
        &= \left(\left(g_1\ast g_2\right)\odot g_3\right)~\widehat{}~(\chi)\\
        &= \left(\left(Uf_1\ast Uf_2\right)\odot Uf_3\right)~\widehat{}~\left(\chi\right).
    \end{align*}
    Hence we have $\left(U\left(\left(f_1\ast f_2\right)\odot f_3\right)\right)~\widehat{}~\left(\chi\right)= \left(\left(Uf_1\ast Uf_2\right)\odot Uf_3\right)~\widehat{}~\left(\chi\right)$ for all $\chi \in \widehat{G}$. From that it follows, $U\left(\left(f_1\ast f_2\right)\odot f_3\right)= \left(Uf_1\ast Uf_2\right)\odot Uf_3.$ Now for $f\in L^2(G,\Hi)$, we have 
    \begin{align*}
        \left(U\left(\tau_g f\right)\right)~\widehat{}~(\chi)& = T\left(\tau_g f\right)\left(\chi\right)\\
        &= \overline{\chi\left(g\right)}Tf(\chi)\\
        &= \overline{\chi(g)} \left(Uf\right)~\widehat{}~(\chi)\\
        &= \left(\tau_g Uf\right)~\widehat{}~(\chi).
    \end{align*}
    Therefore we have, $U\tau_g= \tau_g U$ for all $g\in G$, that is, $U$ commutes with translations. Now for basic tensors $f= F\otimes u\in L^2(G,\Hi), ~~ T(f)= (S\otimes I)(F\otimes u)= SF\otimes u.$ Then define $S_1: L^2(G)\to L^2(G)$ by $S_1(F)= \F^{-1}(S(F))$. Then $Uf= \F^{-1}(Tf)= \F^{-1}(SF\otimes u)= \F^{-1}(SF)\otimes u= S_1(F)\otimes u= (S_1\otimes I)(F\otimes u)$. Thus $U$ agrees with $S_1\otimes I$ on the set of all basic tensors. Now using the fact that closed linear span of the set of all basic tensors is whole of $L^2(G,\Hi)$, we can say that $U= S_1\otimes I.$ Now using theorem \ref{T4.2}, there exists $m\in L^\infty(G)$, such that $(Uf)~\widehat{}~(\chi)= m(\chi)\widehat{f}(\chi)$ for all $f\in L^2(G,\Hi)$. From this we have, $\left(U\left(\left(f_1\ast f_2\right)\odot f_3\right)\right)~\widehat{}~(\chi)= m\left(\chi\right)\left(\left(f_1\ast f_2\right)\odot f_3\right)~\widehat{}~(\chi)$, which implies, $\left\langle \left(Uf_1\right)~\widehat{}~(\chi), \left(Uf_2\right)~\widehat{}~\left(\bar{\chi}\right)\right\rangle \left(Uf_3\right)~\widehat{}~(\chi)= \left\langle \widehat{f_1}(\chi), \widehat{f_2}(\bar{\chi})\right\rangle \left(Uf_3\right)~\widehat{}~(\chi).$ Now for each $\chi\in \widehat{G}$, there exists $h_\chi$ such that $\left(Uh_\chi\right)~\widehat{}~(\chi)\neq 0$. Therefore we have $\left\langle \left(Uf_1\right)~\widehat{}~(\chi), \left(Uf_2\right)~\widehat{}~\left(\bar{\chi}\right)\right\rangle= \left\langle \widehat{f_1}(\chi), \widehat{f_2}\left(\bar{\chi}\right)\right\rangle.$ Now from definition of $U$, we can obtain $\left\langle Tf_1(\chi), T(f_2)(\overline{\chi})\right\rangle=\left\langle \widehat{f_1}(\chi), \widehat{f_2}(\overline{\chi})\right\rangle$. Now using hypothesis (iv), we have $\left\langle T(f_1)(\chi), T(\tilde{f_2})(\chi)\right\rangle= \left\langle \widehat{f_1}(\chi), \widehat{\tilde{f_2}}(\chi)\right\rangle$ and which is true for any $f_1,f_2\in L^2(G,\Hi).$ Therefore we can conclude
    $$\left\langle T(f)(\chi), T(g)(\chi)\right\rangle= \left\langle \widehat{f}(\chi), \widehat{g}(\chi)\right\rangle \quad \text{ for all } f,g\in L^2(G,\Hi).$$
\end{proof}

\section{Cosine Transform on a locally compact abelian group} \label{section5}
\noindent
The cosine transform is a map $\mathcal{C}: L^1(G)\to C_0(COS(G))$ given by the formula 
$$\mathcal{C}(f)(\chi):= \int_G f(x)\chi(x)\, dx,$$
\noindent
where right hand side is defined as Bochner integral and
$$COS(G):=\left\{\chi\in C^b(G): \chi\neq 0, \forall_{x,y\in G} \chi(x)\chi(y)=\frac{\chi(x+y)+\chi(x-y)}{2}\right\}.$$
\noindent
For $G=\mathbb{R}$, we have
$$
COS(G)=\left\{x\mapsto \cos(yx): y\in\mathbb{R}\right\}.
$$
Thus, the notation $COS(G)$ is justified by the fact that its elements are precisely the cosine functions on $\mathbb{R}$.
The cosine convolution $\star_\mathcal{C}: L^1(G)\times L^1(G)\to L^1(G)$ is given by 
$$f\star_\mathcal{C} g(x):= \int_G f(u).\frac{g(x+u)+g(x-u)}{2}\, du.$$
Now we state Krukoswski's result for the characterization of the cosine transform.

\begin{theorem}\cite{krukowski1}
    Let $T: L^1(G) \to C^b(COS(G))$ be a linear and bounded operator. If it satisfies the cosine convolution property 
    $$T(f\star_\mathcal{C}g)=Tf.Tg$$
    for all $f,g\in L^1(G).$ Then there exists an open set $E\subset COS(G)$ and a continuous function $\theta_\mathcal{C}: E\to COS(G)$ such that 
    $$T(f)(\phi)=\begin{cases}
        \mathcal{C}(f)~\circ ~\theta_\mathcal{C}(\phi) & \text{ if } \phi\in E,\\
        \qquad 0 & \text{ if } \phi\in \widehat{G}\setminus E.
    \end{cases}$$
    Furthermore, if we assume that $T$ satisfies the additional property 
    $$T\left(\frac{L_x f+ L_{-x}f}{2}\right)(\phi)=\phi(x) T(f)(\phi)$$
    for every $f\in L^1(G), \phi\in E$ and $x\in G$, then $\theta_\mathcal{C}=\phi.$
\end{theorem}
\noindent
We now define the analogous notion for Hilbert space-valued functions.
We can defined the Hilbert space valued cosine transform $\mathcal{C}_H: L^1(\R,\Hi)\to C_0(COS_H(G))$ by
$$ \mathcal{C}_H(f)(\chi)= \int_G f(x)\chi(x)\, dx,$$
where 
$$
    COS_H(G):= \Bigg\{\chi \in C^b(G): \chi\neq 0, \text{ for all } x,y \in G, \chi(x)\chi(y)= \frac{\chi(x+y)+\chi(x-y)}{2}\Bigg\}.
$$
For two vector valued function, the cosine convolution  $\star: L^1(G,\Hi)\times L^1(G,\Hi)\to L^1(G,\C)$ can be defined by 
$$f\star g(x):= \int_G \left\langle f(u)~,~ \frac{g(x+u)+g(x-u)}{2} \right\rangle\, du$$
Now it can be easily checked that for $f,g\in L^1(G,\Hi)$ we have, $\mathcal{C}_H(f\star g)(\chi)= \langle \mathcal{C}_Hf(\chi), \mathcal{C}_Hg(\bar{\chi})\rangle.$ Then for a scalar valued and a vector valued function, the cosine convolution  $\odot: L^1(G,\C)\times L^1(G,\Hi)\to L^1(G,\Hi)$ can be defined by
$$f\odot g(x)= \int_G f(u)~ \frac{g(x+u)+g(x-u)}{2})\, du.$$
Then it is easy to verify that $\mathcal{C}_H(f\odot g)(\chi)= \mathcal{C}_Hf(\chi)\mathcal{C}_Hg(\chi).$ From these observations the following proposition follows.

\begin{propn}
    Let $f,g,h\in L^1(G,\Hi)$ then $(f\star g)\odot h \in L^1(G,\Hi)$ and
    $$\mathcal{C}_H((f\star g)\odot h)(\chi)= \langle \mathcal{C}_Hf(\chi), \mathcal{C}_Hg(\bar{\chi})\rangle \mathcal{C}_Hh(\chi).$$
\end{propn}
Now in the following proposition, we will observe how Cosine transform behaves on translations.

\begin{propn}
    Let $f\in L^1(G,\Hi)$ and $x \in G$, then $\mathcal{C}_H\left(\frac{L_x f+ L_{-x} f}{2}\right)(\chi)= \chi(x)\mathcal{C}_H(f)(\chi)$ for all $\chi \in COS(G).$
\end{propn}

\begin{proof}
    Let $f\in L^1(G,\Hi), x\in G$ and $\chi\in COS(G)$. Now we calculate the following:
    \begin{align*}
        &\mathcal{C}_H\left(\frac{L_x f+ L_{-x} f}{2}\right)(\chi) = \int \frac{L_xf+L_{-x}f}{2}(y)\chi(y)\, dy= \int \frac{f(y-x)+f(y+x)}{2} \chi(y)\, dy\\
        &= \frac{1}{2}\int f(y-x)\chi(y)\, dy+ \frac{1}{2}\int f(y+x)\chi(y)\, dy\\
        &= \frac{1}{2}\chi(x)\int f(y)\chi(y)\, dy+ \frac{1}{2}\chi(-x)\int f(y)\chi(y)\, dy\\
        &= \frac{1}{2}\chi(x) \mathcal{C}_Hf(\chi)+\frac{1}{2}\chi(-x)\mathcal{C}_Hf(\chi)\\
        &= \frac{\chi(x)+\chi(-x)}{2}\mathcal{C}_Hf(\chi)\\
        &= \chi(x)\mathcal{C}_Hf(\chi).
    \end{align*}
    In the last step we have used $\chi\in COS(G)$, so $\frac{\chi(x)+\chi(-x)}{2}= \chi(x).$
\end{proof}
Now we state our result.
\begin{theorem}
    Let $T: L^1(G,\Hi) \to C^b(COS(G),\Hi)$ be a linear and bounded operator satisfying the following:
    \begin{itemize}
        \item [(i)] $T\left(\left(f\star g\right)\odot h\right)(\chi)= \langle Tf(\chi), Tg(\overline{\chi})\rangle Th(\chi) $ for all $\chi \in COS(G)$, $f,g,h\in L^1(G,\Hi)$\\
        \item [(ii)] $T\left(\frac{L_x f+ L_{-x} f}{2}\right)(\phi)= \phi(x).T(f)(\phi)~~$ for every $f\in L^1(G,\Hi), \phi\in COS(G)$ and $x\in G.$
    \end{itemize}
    Then there exists an open set $E\subset COS(G)$ such that 
    $$\langle Tf(\phi), Tg(\bar{\phi})\rangle= \chi_E(\phi)\langle \mathcal{C}_Hf(\phi), \mathcal{C}_Hg(\bar{\phi})\rangle$$ 
\end{theorem}

\begin{proof}
    Let us define $E:= \{\phi \in COS(G): T(f)(\phi)\neq 0 \text{ for some } f\in L^1(G,\Hi) \}.$ It is easy to see that $E$ is open. Now if $\phi \notin E$, then we are done. But if $\phi \in E$, then there exists $g\in L^1(G,\Hi)$ such that $T(g)(\phi)\neq 0$. Therefore $\langle Tg(\phi), e_k \rangle\neq 0$ for some $k$. Now without loss of generality we can take $\langle Tg(\phi), e_k\rangle=1$. Now we calculate the following:
    \begin{align*}
        &\left\langle Tf_1\left(\phi\right), Tf_2\left(\bar{\phi}\right)\right\rangle\\
        &= \left\langle Tf_1\left(\phi\right), Tf_2\left(\bar{\phi}\right)\right\rangle \left\langle Tg\left(\phi\right), e_k\right\rangle\\
        &= \left\langle \left\langle Tf_1\left(\phi\right), Tf_2\left(\bar{\phi}\right)\right\rangle Tg\left(\phi\right), e_k \right\rangle\\
        &= \left\langle T\left(\left( f_1\ast f_2\right)\odot g\right)(\phi), e_k\right\rangle\\
        &= \left\langle T_\phi\left(\left( f_1\ast f_2\right) \odot g\right), e_k\right\rangle\\
        &= \left\langle T_\phi\left( \int\int \left\langle f_1(z), \frac{f_2(y+z)+f_2(y-z)}{2}\right\rangle \frac{L_y g+ L_{-y} g}{2}\, dy\, dz\right) , e_k \right\rangle\\
        &= \left\langle  \int\int \left\langle f_1(z), \frac{f_2(y+z)+f_2(y-z)}{2}\right\rangle T_\phi\left(\frac{L_y g+ L_{-y} g}{2}\right)\, dy\, dz , e_k \right\rangle\\
        &= \left\langle  \int\int \left\langle f_1(z), \frac{f_2(y+z)+f_2(y-z)}{2}\right\rangle \phi(y)Tg(\phi)\, dy\, dz , e_k \right\rangle\\
        &= \left\langle  \int\int \left\langle f_1(z), \frac{f_2(y+z)}{2}\right\rangle \phi(y)Tg(\phi)\, dy\, dz + \int\int \left\langle f_1(z), \frac{f_2(y-z)}{2}\right\rangle \phi(y)Tg(\phi)\, dy\, dz , e_k \right\rangle\\
         &= \left\langle  \int\int \left\langle f_1(z), \frac{f_2(y)}{2}\right\rangle \phi(y-z)Tg(\phi)\, dy\, dz + \int\int \left\langle f_1(z), \frac{f_2(y)}{2}\right\rangle \phi(y+z)Tg(\phi)\, dy\, dz , e_k \right\rangle\\
         &= \left\langle \left\langle \int f_1(z)\overline{\phi(z)}\, dz, \int \frac{f_2(y)}{2}\overline{\phi(y)}\, dy\right\rangle Tg(\phi)+ 
         \left\langle \int f_1(z)\phi(z)\, dz, \int \frac{f_2(y)}{2}\overline{\phi(y)}\, dy\right\rangle Tg(\phi), e_k \right\rangle\\
         &=\left\langle \left\langle \mathcal{C}_H(f_1)(\bar{\phi}),\frac{1}{2}\mathcal{C}_H(f_2)(\bar{\phi})\right\rangle Tg(\phi)+ 
         \left\langle \mathcal{C}_H(f_1)(\phi), \frac{1}{2}\mathcal{C}_H(f_2)(\bar{\phi})\right\rangle Tg(\phi), e_k \right\rangle\\&= \left\langle \frac{1}{2}\mathcal{C}_H(f_1)(\bar{\phi}+\phi), \mathcal{C}_H(f_2)(\bar{\phi})\right\rangle\left\langle Tg(\phi), e_k\right\rangle\\
         &= \left\langle \mathcal{C}_H(f_1)(\phi), \mathcal{C}_H(f_2)(\bar{\phi})\right\rangle.
    \end{align*}
    Here in the last step, we have used $\frac{\phi+\bar{\phi}}{2}= \phi$, as $\phi \in COS(G)$.  Therefore we obtain, $$\left\langle Tf_1(\phi), Tf_2(\bar{\phi})\right\rangle= \chi_E(\phi)\left\langle \mathcal{C}_Hf_1(\phi), \mathcal{C}_Hf_2(\bar{\phi})\right\rangle.$$
\end{proof}

\section{Laplace Transform on a locally compact abelian group} \label{section6}

In this section, we will see that, like the Fourier transform, the Laplace transform also possesses a similar characterization. First we discuss the result of Krukowski for scalar valued functions.  The Laplace transform is a map $\mathcal{L} : L^1\left(\R_+\right)\to C_0\left(\R_+\right)$ given by 
$$ \mathcal{L}(f)(y):= \int_0^\infty e^{-yx}f(x)\, dx.$$
The Laplace convolution $\star_{\mathcal{L}}: L^1\left(\R_+\right)\times L^1\left(\R_+\right)\to L^1\left(\R_+\right)$ is given by
$$f\star_{\mathcal{L}}g(x):= \int_0^x f(u)g(x-u)\, du.$$
Then it is easy to observe that $\mathcal{L}(f\star_{\mathcal{L}}g)= \mathcal{L}(f)\mathcal{L}(g),~~$  for every $f,g \in L^1(\R_+).$ Using this property Krukowski \cite{krukowski1} characterized the Laplace transform. We will state his result as follows.

\begin{theorem}
    Let $ T : L^1\left(\mathbb{R}_+\right) \longrightarrow C^b\left(\mathbb{R}_+\right) $
be a linear and bounded operator. If it satisfies the Laplace convolution property
   $$ T(f \star_\mathcal{L} g) = T(f) \cdot T(g) $$
for every \( f, g \in L^1\left(\mathbb{R}_+\right) \), then there exists an open set \( E \subset \mathbb{R}_+ \) and a continuous function  $\theta_\mathcal{L} : E \longrightarrow \mathbb{R}_+ \cup \{0\} $
such that
\[
T(f)(y) = 
\begin{cases}
    \mathcal{L}(f)(\theta_{\mathcal{L}}(y)) & \text{if } y \in E, \\
    0 & \text{if } y \in \mathbb{R}_+ \setminus E.
\end{cases}
\]
Furthermore, if we assume that \( T \) satisfies the additional property that there exists \( t_* \in \mathbb{R}_+ \) such that
\begin{equation} \label{eq:shift}
    T(L_{t_*}f)(y) = e^{-y t_*} T(f)(y)
\end{equation}
for every \( f \in L^1\left(\mathbb{R}_+\right) \) and \( y \in E \), where
\[
L_{t_*}f(x) = 
\begin{cases}
    f(x - t_*) & \text{if } x > t_*, \\
    0 & \text{if } 0 < x \leq t_*,
\end{cases}
\]
then \( \theta_{\mathcal{L}}(y) = y. \)
\end{theorem}

Additionally, if we assume $T(L_{t}f)(y) = e^{-y t} T(f)(y)$ for all $t>0,$ then the above proof can be written as follows:

\noindent
Let $E:=\{y\in \R_+: T(f)(y)\neq 0 \text{ for some } f\in L^1(\R_+) \}$. Now if $y\in E$, consider nonzero linear functional $T_y: L^1(\R_+)\to \R$ by $T_y(f)=T(f)(y)$. We choose $g\in L^1(\R_+)$ such that $T_y(g)=1$. Now we calculate
\begin{align*}
    &T(f)(y) = T(f)(y)T(g)(y)= T(f\ast_\mathcal{L}g)(y)= T_y\left(\int_0^\infty f(u) L_u g\, du\right)\\
    &= \int_0^\infty f(u)T_y\left(L_u g\right)\, du= \int_0^\infty f(u)e^{-uy}T(g)(y)\, du= \int_0^\infty f(u)e^{-yu}\, du= \mathcal{L}(f)(y).
\end{align*}
Next, we shall discuss the Laplace transform for Hilbert space valued function. 
The Hilbert space valued Laplace transform $\mathcal{L}: L^1\left(\R^+,\Hi\right)\to C_0\left(\R^+,\Hi\right)$ is defined by 
$$\mathcal{L}(f)(t)= \int_0^\infty e^{-tu}f(u)\, du,$$ where the right hand side is a Bochner integral. 
Now the Laplace convolution of $\Hi$- valued function can be defined as follows.
\begin{equation}
    f\ast g(u):= \int_0^u \langle f(s), g(u-s)\rangle\, ds, ~~ f,g \in L^1(\R^+, \Hi). 
\end{equation}
It is easy to observe that $f\ast g\in L^1(\R^+,\C)$.
Then the Laplace convolution of scalar valued function with a vector valued function can be defined as follows:
\begin{equation}
    f\odot g(u):= \int_0^u f(s)g(u-s)\, ds, ~~ f\in L^1(\R^+, \C), g\in L^1(\R^+, \Hi).
\end{equation}
Now we can check that $f\odot g \in L^1(\R^+,\Hi)$. So, finally we have $(f\ast g) \odot h \in L^1(\R^+,\Hi)$ whenever $f,g,h\in L^1(\R^+,\Hi).$ Then we have the following convolution theorem.
\begin{theorem}
    Let $f,g,h\in L^1(\R^+,\Hi)$, then we have 
    $$\mathcal{L}((f\ast g)\odot h)(t)=\left\langle \mathcal{L} f(t), \mathcal{L} g(t)\right\rangle \mathcal{L} h(t),\quad \text{ for all } t>0.$$
\end{theorem}
\noindent
We have the action of Laplace transform on the translation operator as,

$$\mathcal{L}(L_s f)(t)= e^{-st}\mathcal{L} f(t) \text{ for all } f \in L^1(\R^+,\Hi).$$
\noindent
Now we characterize Laplace transform in this Hilbert space valued set up.

\begin{theorem}
    Let $T: L^1(\R^+, \Hi)\longrightarrow C^b(\R^+, \Hi)$ be a linear and bounded map satisfying the followings:
    \begin{itemize}
        \item [(i)] $T((f\ast g)\odot h)(t)= \langle T f(t), T g(t)\rangle T h(t)$\quad for all $f,g,h \in  L^1(\R^+, \Hi).$
        \item [(ii)] $T(L_s f)(t)= e^{-st} T f(t)$ \quad for all $f\in L^1(\R^+, \Hi), s>0.$
        \item [(iii)] If $\langle f(x),e_k\rangle =0$ for a.e. $x\in \R$,  for some $k$ then $\langle Tf(\xi),e_k\rangle =0$ \quad for all $\xi \in \R^+,$
    \end{itemize}
    then there exist $E\subset \R$, such that 
    
    $$ \left\langle Tf(y), Tg(y)\right\rangle= 
    \begin{cases}
        \left\langle \mathcal{L}f(y), \mathcal{L}g(y)\right\rangle & \text{ if } y\in E,\\
        0 & \text{ if } y\notin E,
    \end{cases} $$  for all $f, g\in L^1(\R^+, \Hi),~~ y\in \R^+.$
\end{theorem}

\begin{proof}
    Let $E:= \{ y\in \R^+: \text{ there exists } f \in L^1(\R^+, \Hi) \text{ such that } T(f)(y)\neq 0\}$. Then it is easy to check that $E$ is open. If $y \notin E$, then we are done. If $y\in E$. Then there exists $h\in L^1(\R^+, \Hi)$ such that $T(h)(y)\neq 0$. Then $\langle T(h)(y), e_k \rangle \neq 0$ for some $k \in \N.$ Now we define $T_y: L^1(\R^+, \Hi)\to \R$ by $T_y(f):= \langle T(f)(y), e_k\rangle.$ Then $T_y$ is nonzero linear functional, that is, $T_y\in \left(L^1(\R^+, \Hi)\right)^\ast$. Now by $L^1-L^\infty$ duality , there exist $\chi_y \in L^\infty (\R^+, \Hi) $ , nonzero such that $T_y(f)= \int_0^\infty \langle f(x), \chi_y(x)\rangle \, dx$ for all $f\in L^1(\R^+, \Hi).$ Now we claim that $\chi_y \perp e_i$ for a.e. $t\in \R$ and $i\neq k.$ If not there exists $j \neq k$ and a measurable set $A$ such that $0< \mu (A)\leq \infty$ and $\langle e_j, \chi_y(t)\rangle\neq 0$, $t\in A.$ Now we define 
    $$\chi_j(t)=\begin{cases}
        \langle e_j, \chi_y(t)\rangle, &  t\in A\\
        0, & o/w.
    \end{cases}$$
    Then $\chi_j \in L^1(\R)$ and non vanishing on a set of positive measure. Then there exist $B\subset \R$ with $0<\mu(B)<\infty$ such that $\int_B \chi_j(t)\neq 0.$ Now let $w(t)= \mathbb{1}_B(t)e_j$, then $w\in L^1(\R,\Hi).$ Now $T_y(w)= \langle T(w)(y), e_k\rangle=0$, by orthogonality hypothesis. But in other way, $T_y(w)= \int_0^\infty \langle w(t), \chi_y(t)\rangle\, dt= \int_B \langle e_j, \chi_y(t)\rangle \, dt\neq 0.$ which is a contradiction. Hence $\chi_y\perp e_i$ for a.e. $t\in \R$ for all $i\neq k.$ So we can assume that $\chi_y(t)= \overline{\psi_y(t)}e_k$, for some $\psi_y\in L^\infty(\R).$ Now we can write $T_y(f)= \int_0^\infty \langle f(t), \chi_y(t)\rangle \, dt= \int_0^\infty \langle f(t), \overline{\psi_y(t)}e_k\rangle\, dt= \left\langle \int_0^\infty \psi_y(t)f(t)\, dt, e_k\right\rangle$. Now for $f,g,h \in L^1(\R, \Hi)$ we have :
    \begin{align*}
        T_y\left(\left(f\ast g\right)\odot h\right) &= \left\langle \int_0^\infty \psi_y(t)\left(f\ast g\right)\odot h(t)\, dt, e_k\right\rangle\\
        &= \left\langle \int_0^\infty \psi_y(t)\int_0^\infty \left( f\ast g\right)(s)h\left(t-s\right)\, ds\, dt, e_k\right\rangle\\
        &= \left\langle \int_0^\infty f\ast g(s)\int_0^\infty \psi_y(t) L_s h(t)\, dt \, ds, e_k\right\rangle\\
        &= \int_0^\infty f\ast g(s) \left\langle \int_0^\infty \psi_y(t) L_s h(t)\, dt , e_k\right\rangle \, ds\\
        &= \int_0^\infty  f\ast g(s) T_y(L_s h)\, ds\\
        &= \int_0^\infty  f\ast g(s) \left\langle T(L_s h)(y), e_k\right\rangle \, ds\\
        &= \int_0^\infty f\ast g (s) \left\langle e^{-sy} T(h)(y), e_k\right\rangle \, ds\\
        &= \left\langle T h(y), e_k\right\rangle \int_0^\infty f\ast g(s) e^{-sy}\, ds\\
        &= \left\langle T h(y), e_k\right\rangle \mathcal{L}(f\ast g)(y)\\
        &= \left\langle T h(y), e_k\right\rangle \left\langle \mathcal{L}f(y), \mathcal{L}g(y)\right\rangle.
    \end{align*}
    Again on the other hand we have,
    \begin{align*}
        T_y\left(\left(f\ast g\right)\odot h\right) &= \left \langle  T\left(\left(f\ast g\right)\odot h\right)(y), e_k\right\rangle\\
        &= \left\langle \left\langle T f(y), T g(y)\right\rangle T h(y), e_k \right\rangle\\
        &= \left\langle Tf(y), Tg(y)\right\rangle \left\langle Th(y), e_k \right\rangle.
    \end{align*}
    As $\langle Th(y), e_k \rangle \neq 0$, we obtain $ \left\langle Tf(y), Tg(y)\right\rangle= \left\langle \mathcal{L}f(y), \mathcal{L}g(y)\right\rangle, $  for all $f, g\in L^1(\R^+, \Hi),~~ y\in \R^+.$
    This completes the proof.
\end{proof}

\noindent{\bf Further Remark.}
Note that the group Fourier transform of Hilbert space-valued functions on the Heisenberg group was studied in \cite{swain2010p}.
Let $f:\He \to \Hi$ be an integrable function, then the group Fourier transform of $f$ is defined by
$$\widehat{f}(\lambda)=\int_{\He}f(w,s)\pi_\lambda(w,s)\, dw\, ds, \qquad \lambda\in \R^\ast, $$
where the right hand side is defined as a Bochner integral.

Now let $\varphi\in L^1(\He,\C), h\in L^1(\He,\Hi)$, then the $\odot$- convolution is defined by
$$\varphi\odot h(z,t)=\int_{\He} \varphi(w,s)h((z,t)(-w,-s))\, dw\, ds.$$
Then one can easily check that $\varphi\odot h\in L^1(\He, \Hi).$
Here, we obtain a characterization of the group Fourier transform on the Heisenberg group when $\Hi=\C$. It would be interesting to know if the above characterizations results hold for the group Fourier transform of Hilbert space-valued functions on the Heisenberg group. 

\smallskip

\noindent{\bf Acknowledgements.} 

\noindent
This work is a part of first author's Ph.D. thesis. He thanks the National Institute of Science Education and Research Bhubaneswar, India, for providing an excellent research facility. RM acknowledges the partial support provided by the research grants  (ANRF/ARGM/2025/000750/MTR).

\smallskip
\noindent{\bf Funding.}  Not applicable.

\smallskip
\noindent{\bf Data Availability Statement.} Not applicable.

\smallskip
\noindent{\bf Code Availability Statement.} Not applicable.

\smallskip
\noindent{\bf Declarations.}

\noindent{\bf Conflict of interest.} The authors declare that they have no conflict of interest.
\smallskip


\nocite{krukowski2}
\nocite{Kumar}
\nocite{Lavanya_Weyl}
\nocite{Lavanya_Heisen}
\medskip

\bibliography{ref1.bib}
\bibliographystyle{amsplain}

\end{document}